\documentclass[12pt]{article}
\usepackage{amsmath}
\usepackage{amsthm}
\usepackage{amsfonts}
\usepackage{amssymb}
\usepackage{hyperref}
\usepackage{mathrsfs}

\theoremstyle{definition}

\newcommand{\s}{{\tt switch}}
\newcommand{\m}{{\tt match}}
\newcommand{\D}{{\mathcal D}}

\begin{document}
\title{Dominance in the Monty Hall Problem}
\author{Alexander V. Gnedin\thanks{Postal address:
 Department of Mathematics, Utrecht University,
 Postbus 80010, 3508 TA Utrecht, The Netherlands. E-mail address: A.V.Gnedin@uu.nl}}
\date{}
\maketitle

\begin{abstract}
\noindent
Elementary decision-theoretic analysis of the Monty Hall dilemma  shows that the problem has dominance.
This makes possible 
to discard nonswitching strategies, without making any assumptions on the prior distribution of factors out of control of the decision maker.
A path to the Bayesian and the minimax decision-making environments  is then straightforward.

\end{abstract}
\noindent
Suppose there is a set of `doors' $\D$ with at least three elements.
One door $\theta\in\D$ is a winning door. You will be asked to choose a door $x\in\D$, and then all doors except one door
$y\in\D\setminus \{x\}$ will be revealed as not winning. Then you will be asked to guess if  $x=\theta$ (action \m~ with door $x$) or
$x\neq \theta$ (action \s~ to door $y$). You win if your guess is correct.

The case of three doors is widely known as the Monty Hall problem.
The observation of this note is that the problem has {\it dominance}, 
which makes sense of the question `to switch or not to switch' even in the situations
when no prior distributions are assigned to the parameters out of control of the decision maker. 
The idea is quite simple: if  the guessing strategy chooses $x$ and in some situation plays $\m$, so does not
switch to $y$, then  
another strategy which chooses $y$ and plays $\s$ all the time is at least as good as the first whichever $\theta\in \D$.
This simplistic view is basically right, 
modulo subtleties involved in the formal definition of  
`strategy' and
`situation' 

A {\it strategy} is a pair $(x,a_x)$ where $x\in\D$ and $a_x$ is a function of $y\in \D\setminus\{x\}$
with values $a_x(y)\in\{\m,\s\}$. 
For constant functions
we call $(x, \m)$ and $(x, \s)$ {\it single-action} strategies, and call $(x,\s)$  {\it always-switching} strategy.

The decision-making environment is specified
by two parameters.
 One parameter is  the winning door $\theta\in \D$. Another parameter
is the door $y$ to which switching can occur, which is a
 function $d_\theta$  of $x\in \D$ , that must be {\it admissible} in the sense that   $d_\theta(x)=\theta$ off the diagonal $\theta=x$,
and $d_\theta(x)\neq \theta$ on the diagonal $\theta=x$.

Let $W(\theta,x,{\tt match}):=1(x=\theta)$ and $W(\theta,x,{\tt switch}):=1(x\neq\theta),$
where $1(\cdots)$ denotes indicators. 
A strategy $(x, a_x)$ 
is evaluated by the win-or-nothing payoff function $W(\theta, x,a_x(d_\theta(x)))$.
Note that the payoff does not depend on $y=d(\theta,x)$ explicitly.
The payoff of each single-action strategy depends on $y$ 
neither explicitly nor implicitly.

\vskip0.2cm
{\bf Theorem} {\it The class of always-switching strategies is weakly dominant, that is
for every strategy $(x,a_x)$ there exists an always-switching strategy $(x',\s)$ such that 
$$W(\theta,x',\s)\geq W(\theta,x,a_x(d_\theta(x)))$$ 
for all $\theta\in\D$ and all admissible functions $d_\theta$.
}

\begin{proof}
 If $a_x$ is not always-switching, there exists $x'\neq x$ with
 $a_x(x')={\tt match}$.  
Consider the strategy $(x',\s)$. For $\theta\neq x'$ we have 
$$W(\theta,x,a_x(d_\theta(x)))\leq 1= W(\theta,x',\s).$$
For  $\theta= x'$ we have 
$W(x',x',\s)=0$,  
but then    also
$$W(\theta,x,a_x(d_\theta(x)))=  W(x',x,a_x(x'))= W(x',x,\m)=0,$$
since by the admissibility $d_\theta(x)=x'$. 
\end{proof}

In the Bayesian setting of the decision problem the winning door  $\Theta$ and the door $Y$   to which switching is offered 
are random variables with given probability
distributions, with distribution of $Y$ defined conditionally on the decision variable $x$ and $\Theta$, with account 
of the rule $Y=\Theta$ in the event  $\Theta\neq x$.
Assuming $\D$ finite, the distributions can be specified by a probability mass function $p_\theta$ and
transition probabilities $q_{\theta,y}$ for $\theta\in\D, y\in\D\setminus\{\theta\}$.
A Bayesian strategy is the always-switching strategy $(\theta^*, \s)$, where $\theta^*$ minimizes $p_\theta$.
The optimality readily follows by first noting that by the dominance the class of always-switching strategies is optimal,
then noting that the winning probability with $(x,\s)$ is $1-p_x$, {\it i.e.} the probability of the event $\Theta\neq x$ 
that the winning door is not the one you have chosen before facing the `to switch or not to switch' dilemma.

Let $n$ be the cardinality of $\D$ ($n\geq 3$). Since $\max_{x\in\D}(1-p_x)\geq (n-1)/n$  it is clear that in the worst-case  
the distribution of $\Theta$ is  uniform  on $\D$, and the minimax probability of winning the prize is $v:=(n-1)/n$.
If $\Theta$ is uniform, then every always-switching strategy yields $v$, no matter what is the distribution of $Y$.
On the other hand, if the decision-maker applies the mixed strategy $(X,\s)$ with $X$ uniformly distributed, then 
the probability to win is $v$ for arbitrary $(\Theta,Y)$ (independent of $X$).
It follows that in the zero-sum setting of the decision problem, every saddle-point solution has the form
$(X,\s)$-versus-$(\Theta,Y)$,  
with $X$ and $\Theta$ uniformly distributed and the law of $Y$ 
determined by arbitrary  transition probabilities. The value of this game is then $v=(n-1)/n$.

\vskip0.5cm
\noindent
{\small The Monty Hall problem (MHP) and its numerous variants, all involving probabilities, are summarized on the Wikipedia MHP page.
We refer to \cite{Gill} for a recent critical analysis of the textbook solutions and the history of the MHP.

\vskip0.3cm
\noindent
{\bf Acknoweledgement.} The author is indebted to Richard Gill for motivation. 

\end{document}